 \documentclass[12pt]{amsart}
 \usepackage{amssymb}
 \usepackage{amsmath,amsxtra}
\usepackage{amsfonts}

\numberwithin{equation}{section}

 \newtheorem{thm}{Theorem}[section]
 \newtheorem{lemma}[thm]{Lemma}
 \newtheorem{prop}[thm]{Proposition}
 \newtheorem{cor}[thm]{Corollary}

 \theoremstyle{definition}

\newtheorem{problem}{Problem}

\theoremstyle{remark}
\newtheorem{remark}[thm]{Remark}

\newcommand{\Sc}{\mathcal{S}}

\newcommand{ \Rb}{\mathbb{R}}

 \newcommand{ \Cbb}{\mathbb{C}}

 \newcommand{ \Nbb}{\mathbb{N}}
 \newcommand{ \Zbb}{\mathbb{Z}}

 \newcommand{\al}{\alpha}
 \newcommand{\la}{\langle}
 \newcommand{\ra}{\rangle}

 \newcommand{\ep}{\varepsilon}

\newcommand{\lamd}{\lambda}

 \newcommand{\LO}{\mathcal{LO}}
\newcommand{\p}{\,|\,}

\begin{document}

 \title[Fixed point properties on unbounded sets]{Fixed point properties for semigroups of nonlinear mappings on unbounded sets}

\author[A. T.-M. Lau]{Anthony To-Ming Lau \dag}
\address{\dag \; Department of Mathematical and Statistical sciences\\
           University of Alberta\\
           Edmonton, Alberta\\
           T6G 2G1 Canada}
\email{tlau@math.ualberta.ca}
\thanks{\dag  \; Supported by NSERC Grant MS100}

\author[Y. Zhang]{ Yong Zhang \ddag}
\address{\ddag \; Department of Mathematics\\
           University of Manitoba\\
           Winnipeg, Manitoba\\
           R3T 2N2 Canada}
\email{zhangy@cc.umanitoba.ca}
\thanks{\ddag \; Supported by NSERC Grant 238949-2005}

\date{November 2014}



\begin{abstract}
A well-known result of W. Ray asserts that if $C$ is an unbounded convex subset of a Hilbert space, then there is a nonexpansive mapping $T$: $C\to C$ that has no fixed point. In this paper we establish some common fixed point properties for a semitopological semigroup $S$ of nonexpansive mappings acting on a closed convex subset $C$ of a Hilbert space, assuming that there is a point $c\in C$ with a bounded orbit and assuming that certain subspace of $C_b(S)$ has a left invariant mean. Left invariant mean (or amenability) is an important notion in harmonic analysis of semigroups and groups introduced by von Neumann in 1929 \cite{Neu} and formalized by Day in 1957 \cite{Day}. In our investigation we use the notion of common attractive points introduced recently by S. Atsushiba and W. Takahashi.
\end{abstract}

\maketitle

\section{Introduction}

Let $E$ be a Banach space and $C$ be a nonempty bounded closed convex subset of $E$. The set $C$ is said to have the fixed point property (abbreviated as fpp)  if every nonexpansive mapping $T$: $C\to C$ has a fixed point, where $T$ being nonexpansive means $\|T(x) - T(y)\| \leq \|x-y\|$ for all $x,y\in C$. The space $E$ is said to have the fpp if every bounded closed convex set of $E$ has the fpp.

A result of Browder \cite{Br} asserts that if a Banach space $E$ is uniformly convex, then $E$ has the fpp. As shown by Aspach \cite{Alspach} (see also \cite[Example~11.2]{G-K}), there is a weakly compact convex subset of $L^1[0,1]$ on which an isometry does not have a fixed point. It is also well-known that a weak* compact convex subset of $\ell^1(\Zbb)$ has the fpp. However, $\ell^1(\Zbb)$ does not have the fpp for bounded closed convex sets \cite{G-K}. In a recent remarkable paper of Lin \cite{Lin}, it was shown that $\ell^1(\Zbb)$ can be renormed to have the fpp. This answers in negative a long-standing open question of whether every Banach space with the fpp is necessarily reflexive. It was proved by  B. Maurey in \cite{Maurey} that every nonempty weakly compact convex subset of the sequence space $c_0$ has the fpp for nonexpansive mappings. In the beautiful paper \cite{Benavides}, T. D. Benavides proved that for every unbounded subset $C$ in $c_0$ there is a nonexpansive mapping $T$ on $C$ which is fixed point free.


Let $S$ be a semitopological semigroup, that is, a semigroup with a Hausdorff topology such that for each $t\in S$, the mapping $s\mapsto t\cdot s$ and $s\mapsto s\cdot t$ from $S$ into $S$ are continuous. Let $C$ be a subset of a 
Banach space $E$. We say that $\Sc = \{T_s:\; s\in S\}$ is a \emph{representation} of $S$ on $C$ if for each $s\in S$, $T_s$ is a mapping from $C$ into $C$ and $T_{st}(x) = T_s(T_tx)$ ($s,t\in S$, $x\in C$). Sometimes we simply use $sx$ to denote $T_s(x)$ if there is no confusion in the context. The representation is called \emph{separately} or, respectively, \emph{jointly continuous} if the mapping $(s,x)\mapsto T_s(x)$ from $S\times C$ to $C$ is separately or jointly continuous. We say that a representation $\Sc$ is 
\emph{nonexpansive} if 
$\|T_sx - T_sy\| \leq \|x-y\|$ for all $s\in S$ and all $x,y \in C$. A point $x\in C$ is called a \emph{common fixed point} for (the representation of) $S$ if $T_s(x) = x$ for all $s\in S$. The set of all common fixed points for $S$ in $C$ is called the \emph{fixed point set} of $S$ (in $C$) and is denoted by $F(S)$.

Let $\Sc$ be a jointly continuous representation of $S$ on a closed convex subset $C$ of a Hilbert space $H$. Then, as well-known, $F(S)$ is a closed and convex subset of $C$ if it is not empty \cite{Demarr}. However, $F(S)$ may be empty for a continuous representation of $S$ on an unbounded convex set $C$ of a Hilbert space even if $S$ is a commutative semigroup with a single generator \cite{Ray}.

In the recent paper \cite{A-T}, which was motivated by \cite{T-T}, Atsushiba and Takahashi introduced the concept of common attractive points for a nonexpansive representation $\Sc$ of a semigroup $S$ on a set $C$ in a Hilbert space $H$ (precise definition may be seen in Section~\ref{preliminaries}). They showed that $F(S)\neq \emptyset$ for commutative $S$ if there is a common attractive point for $S$ \cite[Lemma~3.1]{A-T}. They showed further that for commutative semigroups $S$, if $\{T_sc, s\in S\}$ is bounded for some $c\in C\subset H$, then the set $A_C(\Sc)$ of all attractive points of $\Sc$ is not empty. As a consequence, $F(S)\neq \emptyset$ \cite[Theorem~4.1]{A-T}. We note that the assumption that $\{T_sc, s\in S\}$ is bounded for some $c\in C\subset H$ cannot be dropped even when $S$ is commutative. Indeed, by the classical result of  W. Ray in \cite{Ray} as mentioned above, for every unbounded convex subset $C$ of a Hilbert space there is a nonexpansive mapping $T_0$: $C\to C$ without a fixed point in $C$. In particular, the representation $\{T_n(c): n\in \Nbb\}$ of $(\Nbb, +)$ does not have a common fixed point in $C$. An investigation continuing that of \cite{A-T} may be seen in \cite{T-W-Y}.

As one of the main results in this paper, we show that the result mentioned above of Atsushiba's and Takahashi remains true when $\Sc$ is a continuous representation of a left amenable semitopological semigroup $S$, where $S$ being left amenable means $C_b(S)$ of bounded continuous complex-valued functions on $S$ has a left invariant mean. It also remains true when $S$ is separable and left reversible if the representation is weakly equicontinuous. Here a semitopological semigroup $S$ is \emph{left reversible} if any two closed right ideals of $S$ have non-void intersection, that is, $\overline{sS}\cap \overline{tS} \neq \emptyset$ for any $s,t\in S$, where, for a subset A of a topological space, $\overline{A}$ denotes the closure of $A$. This is the case when $S$ is normal or $C_b(S)$ has a left invariant mean \cite{H-L}. Likewise, $S$ is \emph{right reversible} if any two closed left ideals of $S$ have non-void intersection, that is, $\overline{Ss}\cap \overline{St} \neq \emptyset$ for any $s,t\in S$.

The paper is organized as follows: In Section~\ref{B-spaces} we study the relation between the common attractive point and the common fixed point for a semigroup of nonexpansive mappings on a closed convex subset $C$ of a strictly convex space. In Section~\ref{Hilbert} we establish our main results concerning common fixed points on a closed convex subset of a Hilbert space. In Section~\ref{hybrid} we extend some of our results in Section~\ref{Hilbert} to the class of generalized hybrid mappings introduced recently in \cite{K-T-Y}. In Section~\ref{problems} we post some related open problems.

\section{Preliminaries and notations}\label{preliminaries}

Topologies considered in this paper will be all Hausdorff. Banach spaces are all assumed to be over the complex numbers $\Cbb$. If $E$ is a Banach space (resp. a dual Banach space), the weak topology (resp. weak* topology) of $E$ will be denoted by wk (resp. wk*).

Let $\Sc = \{T_s:\; s\in S\}$ be a representation of a semigroup $S$ on a convex subset $C$ of a Banach space $E$. The representation is called \emph{affine} if $C$ is convex and each $T_s$ ($s\in S$) is an affine mapping, that is,  $T_s(ax+by) = aT_sx + bT_sy$ for all constants $a,b \geq 0$ with $a +b = 1$, $s\in S$ and $x, y\in C$. A point $a\in E$ is an \emph{attractive point} of $\Sc$ 
if $\|a-T_sx\| \leq \|a-x\|$ for all $x\in C$. The set of all attractive points of $\Sc$ for $C$ is denoted by $A_C(\Sc)$.

Given a semitopological semigroup $S$,  let $\ell^\infty(S)$ be the $C^*$-algebra of bounded complex-valued functions on $S$ with the supremum norm and pointwise multiplication. For each $s\in S$ and $f\in \ell^\infty$, denote by $\ell_sf$ and $r_sf$ the left and right translates of $f$ by $s$ respectively, that is, $(\ell_sf)(t) = f(st)$ and $(r_sf)(t) = f(ts)$ ($t\in S$). Let $X$ be a closed subspace of $\ell^\infty(S)$ containing the constant functions and being invariant under left translations. Then a linear functional $m\in X^*$ is called a mean if $\|m\| = m(1) =1$; $m$ is called a left invariant mean, denoted by LIM, if $m(\ell_sf) = m(f)$ for all $s\in S$, $f\in X$. 
Let $C_b(S)$ be the space of all bounded continuous complex-valued functions on $S$. $C_b(S)$ certainly is a closed \-subalgebra of $\ell^\infty(S)$ containing the constant functions and being invariant under translations. Let $LUC(S)$ be the space of left uniformly continuous functions on $S$, that is, all $f\in C_b(S)$ such that the mappings $s\mapsto \ell_s(f)$ from $S$ into $C_b(S)$ are continuous. 
Then $LUC(S)$ is a $C^*$-subalgebra of $C_b(S)$ invariant under translations and contains the constant functions. When $S$ is a topological group, then $LUC(S)$ is precisely the space of bounded right uniformly continuous functions on $S$ as defined in \cite{Greenleaf}.
 The semigroup $S$ is called \emph{left amenable} (respectively extremely left amenable) if $LUC(S)$ has a LIM (respectively a multiplicative LIM). Left amenable semitopological semigroups include all commutative semigroups, all compact groups and all solvable groups. But the free group (or semigroup) on two generators is not left amenable. The theory concerning amenability of semigroups may be found in monographs \cite{Pat} and \cite{Pier}.

Let $AP(S)$ be the space of all $f\in C_b(S)$ such that $\LO(f) = \{\ell_sf:\, s\in S\}$ is relatively compact in the norm topology of $C_b(S)$, and let $WAP(S)$ be the space of all $f\in C_b(S)$ such that $\LO(f)$ is relatively compact in the weak topology of $C_b(S)$. Functions in $AP(S)$ (respectively $WAP(S)$) are called \emph{almost periodic} (respectively \emph{weakly almost periodic}) functions. $AP(S)$ and $WAP(S)$ are closed C*-subalgebras  of $C_b(S)$ invariant under translations and contains the constant functions. In general, the following inclusions hold:
\[ AP(S) \subseteq LUC(S) \subseteq C_b(S) \text{ and }  AP(S) \subseteq WAP(S) \subseteq C_b(S). \]
If $S$ is a discrete semigroup then
\[ AP(S) \subseteq WAP(S)\subseteq LUC(S) =C_b(S) = \ell^\infty(S). \]
If $S$ is a compact semitopological semigroup then 
\[ AP(S) = LUC(S)\subseteq WAP(S) = C_b(S). \]
If $S$ is a compact topological semigroup, that is, the multiplication is jointly continuous, then
\[ AP(S) = WAP(S)= LUC(S) = C_b(S). \]
All inclusions asserted in the above diagrams may be proper (see \cite{B-J-M} for details).

Suppose that $C$ is a subset of a 
Banach space $E$ and that $\Sc = \{T_s:\; s\in S\}$ is a representation of $S$ on $C$. Let $c\in C$ be such that $\{T_sc: s\in S\}$ is bounded. Then each $f\in E^*$ defines an element $f_c \in \ell^\infty(S)$ for which $f_c(s) = \la T_sc, f\ra$ for $s\in S$. If $s\mapsto T_sc$: $S\to C$ is continuous when $C$ is equipped with the weak topology of $E$, then $f_c\in C_b(S)$; if the action of $S$ on $C$ is weakly jointly continuous and $\{T_sc: s\in S\}$ is weakly relatively compact, then $f_c \in LUC(S)$.
 If the action of $S$ on $C$ is weakly separately continuous and weakly equicontinuous continuous and $\{T_sc: s\in S\}$ is weakly relatively compact, then $f_c\in AP(S)$ \cite[Lemma~3.1]{Lau73}. Finally, if the action of $S$ on $C$ is weakly separately continuous and weakly quasi-equicontinuous and $\{T_sc: s\in S\}$ is weakly precompact, then $f_c\in WAP(S)$ \cite[Lemma 3.2]{L-Z1}. Here we recall that a representation $\Sc = \{T_s:\; s\in S\}$ on a Hausdorff space $X$ is quasi-equicontinuous if $\overline{\Sc}^{\,p}$, the closure of $\Sc$ in the product space $X^X$, consists of only continuous mappings. In other words, the representation is quasi-equicontinuous if for any net $(s_i)\subset S$, whenever $T_{s_i}(x) \to T(x)$ for each $x\in X$, $T$ is a continuous mapping from $X$ into $X$. 
 
Now let $X$ be a closed subspace of $\ell^\infty(S)$ containing the constant functions. Let $\Sc = \{T_s:\; s\in S\}$ be a representation of $S$ on $C$ as above. Suppose that $c\in C$ such that $\{T_sc: s\in S\}$ is bounded, and suppose that $f_c\in X$ for each $f\in E^*$. For any mean $\mu\in X^*$ on $X$ we may define $T_\mu c\in E^{**}$ by
\[
\la T_\mu c, f\ra = \mu(f_c).
\]
$T_\mu c$ is clearly well-defined. If $\{T_sc: s\in S\}$ is precompact, then $T_\mu c$ is weak* continuous. So $T_\mu c\in E$ in this case.

\section{Attractive points and common fixed point properties}\label{B-spaces}

Recall that a Banach space $E$ is \emph{strictly convex} if $\|\frac{x+y}{2}\| < 1$ whenever $x,y\in E$, $\|x\| = \|y\| =1$ and $x\neq y$. It is readily seen that for any distinct elements $x,y_1, y_2$ from a strictly convex space with $\|x-y_1\| = \|x-y_2\| = d$ we have $\|x-\frac{y_1+y_2}{2}\| < d$. $E$ is \emph{uniformly convex} if for each $0<\ep\leq 2$ there exists $\delta >0$ such that $\|\frac{x+y}{2}\|< 1-\delta$ whenever $x,y\in E$, $\|x\| = \|y\| = 1$ and $\|x-y\|\geq \ep$. It is known that if $E$ is uniformly convex then it is strictly convex and reflexive. Tipical examples of a uniformly convex space are $L^p$-spaces ($p>1$).

Now suppose $E$ is a strictly convex and reflexive Banach space.
Let $C\neq \emptyset$ be a convex subset of $E$. For any $x\in E$ there is a unique $u\in \overline C$, the norm closure of $C$, such that $\|u-x\|\leq \|c-x\|$ for all $c\in C$. In fact, this is trivial if $x\in \overline C$; if $x\notin \overline C$, let
\[ 
C_x = \{y\in \overline C: \|x-y\|\leq \|x-c\| \text{ for all } c\in C\}.
 \]
Then $C_x \neq \emptyset$ since $C_x = \cap_{\al > d}C_x(\al)$, where $d = \inf\{\|x-c\|: c\in C\}$ and $C_x(\al)=\{c\in \overline C: \|x-c\|\leq \al\}$. $C_x(\al)\neq \emptyset$ is a closed convex bounded subset of $E$ and hence is weakly compact since $E$ is reflexive. Finite intersection property implies that $C_x$ is a non-empty weakly compact convex subset of $E$. Let $y_1,y_2 \in C_x$. then $\frac{y_1+y_2}{2}\in C_x$ and hence $\|x-\frac{y_1+y_2}{2}\|=\|x-y_1\| =\|x-y_2\| =d$. By the strict convexity of $E$ this implies $y_1 = y_2$. Therefore, $C_x$ is a singleton. So the element  $u\in C_x$ is the only element of $\overline C$ that satisfies $\|u-x\|\leq \|c-x\|$ for all $c\in C$. If $C$ is closed then
 we  call  this $u$ the \emph{metric projection} of $x$ in $ C$ and denote it by $P_C(x)$. 
 If $E$ is a Hilbert space then $P_C(x)$ may also be characterized as the unique element $u\in C$ satisfying
\begin{equation}\label{mp}
\text{Re}\la x-u \, | \, u-c \ra \geq 0 \quad c\in C.
\end{equation}

\begin{lemma}\label{As to Fs}
Suppose that $E$ is a strictly convex and reflexive Banach space. Let $C\neq \emptyset$ be a closed convex subset of $E$ and $\Sc$ be a representation of a 
semigroup $S$ on $C$. If $A_C(\Sc) \neq \emptyset$, then $F(S) \neq \emptyset$.
\end{lemma}

\begin{proof}
Take $a\in A_C(\Sc)$. Let $u = P_C(a)$. Since $a$ is attractive,
\[
\|a - T_tu\|\leq \|a-u\| \leq \|a-c\| 
\]
for all $c\in C$ ($t\in S$). Thus $T_tu = P_C(a) = u$ for all $t\in S$, which means $u\in F(S)$.

\end{proof}

\begin{remark}
  If $E$ is a general Banach space, the proof of Lemma~\ref{As to Fs} still works as long as $P_C (a)$ is uniquely defined.
\end{remark}

\begin{remark}
The converse of Lemma~\ref{As to Fs} cannot be true in general. Namely, even $F(\Sc)\neq \emptyset$, it still can happen that $A_C(\Sc) = \emptyset$. For example, Let $E = \ell^p$ ($p\geq 1$), $C = \{x\in E:\, x = (x_i)_{i=1}^\infty, x_1 \geq0\}$. Consider $T((x_i)) = (x_1,x_1, x_2, x_3,\cdots)$. Then $T$ is a nonexpansive mapping on $C$, and $T$ has fixed point $\hat0 = (0,0,0,\cdots)$. But T has no attractive point for $C$. As a consequence, the representation $\{T_n = T^n:\, n\in \Nbb\}$ of $(\Nbb, +)$ has a common fixed point in $C$ but  has no attractive points for $C$.
\end{remark}

However, we have a weaker relation between the existence of a fixed point and the existence of an attractive point for nonexpansive representations of a semigroup as follows.

\begin{prop}\label{equiv 1}
Let $E$ be a reflexive, strictly convex Banach space and $C$ a closed convex subset of $E$. Suppose that $\Sc$ is a representation of a semigroup $S$ on $C$ as nonexpansive self mappings. Then the following statements are equivalent
\begin{enumerate}
\item\label{C0As} There is a closed, $S$-invariant convex subset $C_0$ of $C$ such that $A_{C_0}(\Sc)\neq \emptyset$;
\item\label{Fs} $S$ has a common fixed point in $C$.
\end{enumerate}
\end{prop}

\begin{proof}
Assume (\ref{C0As}) holds. Apply Lemma~\ref{As to Fs} for $C_0$. We then see that (\ref{Fs}) is true.

Suppose (\ref{Fs}) holds. We consider $C_0 = F(S)$. Then $C_0$ is closed and $S$-invariant. The convexity of $C_0$ follows from \cite[Lemma~3.4]{G-K}. As $C_0$ is the fixed point set of $S$, it is obvious that all elements of $E$ are attractive points for $C_0$. So (\ref{C0As}) holds.

\end{proof}

We now consider more general representations.

Let $E$ be a Banach space and $C\subset E$. We call a mapping $T$: $C\to C$ \emph{asymptotically nonexpansive} if for all $x,y\in C$ the following inequality holds.
\[
\limsup_{n\to \infty} \|T^nx - T^ny\| \leq \|x-y\|.
\]
We note that the notion was defined by K. Goebel and W.A. Kirk \cite{G-K 1972} in a slightly different way, where they called $T$ asymptotically nonexpansive if there was a sequence $(k_n)$ of real numbers such that $k_n \to 1$ and 
\[
\|T^nx - T^ny\| \leq k_n \|x-y\| \quad (x,y\in C).
\]
Our definition of asymptotic nonexpansiveness is more general than the above Goebel and Kirk's version of the notion. For example, The mapping $T$: $[0,1] \to [0,1]$ defined by $T(x) = \sqrt{x}$ for $x\neq 0$ and $T(0)=1$ is asymptotically nonexpansive in our definition but not in their definition. It is also different from the notion introduced in \cite{H-L2}. Our asymptotically nonexpansive mappings could be discontinuous.

Suppose $\Sc = \{T_s: s\in S\}$ is a representation of a semigroup $S$ on a set $C$ in a Banach space $E$. We call an element $a\in E$ an \emph{asymptotically attractive point} of $\Sc$ for $C$ if 
\[
\limsup_{n\to \infty} \|a - (T_t)^n(x)\| \leq \|a - x\|
\]
for all $x\in C$ and all $t\in S$. We denote the set of all asymptotically attractive points of $\Sc$ for $C$ by $AA_C(\Sc)$.

Certainly, any attractive point is asymptotically attractive and any nonexpansive mapping is asymptotically nonexpansive. 

\begin{prop}\label{equiv 2}
Let $E$ be a reflexive, strictly convex Banach space and $C$ a closed convex subset of $E$. Suppose that $\Sc$ is a representation of the semigroup $S$ on $C$ as weakly continuous and norm asymptotically nonexpansive self mappings. Then the following statements are equivalent
\begin{enumerate}
\item\label{C0 2} There is a closed, $S$-invariant convex subset $C_0$ of $C$ such that $A_{C_0}(\Sc)\neq \emptyset$;
\item \label{C0 AAs} There is a closed, $S$-invariant convex subset $C_0$ of $C$ such that $AA_{C_0}(\Sc)\neq \emptyset$;
\item\label{Fs 2} $S$ has a common fixed point in $C$.
\end{enumerate}
\end{prop}

\begin{proof}
The implication of (\ref{C0 2})$\Rightarrow$(\ref{C0 AAs}) is trivial. We show (\ref{C0 AAs})$\Rightarrow$(\ref{Fs 2}) and (\ref{Fs 2})$\Rightarrow$(\ref{C0 2}).

Suppose that (\ref{C0 AAs}) holds. Without generality we may assume $AA_C(\Sc)\neq \emptyset$. Take $a\in AA_C(\Sc)$ and let $u=P_C(a)$. Fix $t\in S$. Let $v$ be a weak cluster point of $\{(T_t)^nu)\}$. Then $v\in C$ since as a closed and convex set $C$ is weakly closed. We may assume
\[
v = \text{wk-}\lim_k (T_t)^{n_k} u.
\]
Then
\[
\|a - v\| \leq \limsup_k \|a - (T_t)^{n_k} u\| \leq \limsup_{n\to \infty} \|a - (T_t)^{n} u\|.
\]
Since $a$ is an asymptotically attractive point, the above shows $\|a - v\| \leq \|a - u\|$. By the definition f $u$ we derive
\[
v = u = P_C(a).
\]
This is true for every weak cluster point $v$ of $\{(T_t)^nu)\}$. But $\{(T_t)^nu)\}$ is bounded and hence is a subset of a weakly compact set in $C$ (note $E$ is reflexive).  We then conclude
\[
\text{wk-}\lim (T_t)^n u = u.
\]
Using weak-weak continuity of $T_t$, we finally have
\[
T_t (u) = \text{wk-}\lim _n (T_t)^{n+1} (u) = \text{wk-}\lim (T_t)^n u = u.
\]
The above is true for every $t\in S$. Therefore $u\in F(S)$ and hence (\ref{Fs 2}) is true.

Now assume (\ref{Fs 2}). We prove (\ref{C0 2}) holds. We consider $C_0 = F(S)$. This is clearly a nonempty closed $S$-invariant subset of $C$, and $AA_{C_0}(\Sc) = E$. To complete the proof we only need to show that $C_0$ is convex. To this end it suffices to show $\frac{1}{2}(x+y) \in C_0$ whenever $x,y\in C_0$. 

Let $x,y\in C_0$. We may assume 
\[
d = \|x-y\| >0.
\]
Denote $z = \frac{1}{2}(x + y)$ and let $t\in S$. Suppose that $\tilde z$ is a weak cluster point of $\{(T_t)^n z\}$ and assume
\[
\tilde z = \text{wk-}\lim_k (T_t)^{n_k} (z).
\]
Then 
\[
\|\tilde z - x\| \leq \limsup\|(T_t)^n (z) -x\| = \limsup\|(T_t)^n (z) -(T_t)^n(x)\|.
\]
Since the representation of $S$ on $C$ is asymptotically nonexpansive, the above leads to $\|\tilde z -x\| \leq \|z - x\| = \frac{1}{2}d$. Similarly, $\|\tilde z -y\| \leq \frac{1}{2}d$. So 
\[
d = \|x-y\|\leq \|\tilde z - x\| + \|\tilde z - y\| \leq d.
\]
Thus
\[
\|\tilde z - x\| = \|\tilde z - y\| = \frac{1}{2} d.
\]
Let $c = \frac{1}{2}(z+\tilde z)$. Then
\[
\|c - x\| \leq \frac{1}{2}(\|z-x\| +\|\tilde z -x\|) = \frac{1}{2}d, \text{ and similarly } \|c- y\| \leq \frac{1}{2}d.
\]
These show that $\|c-x\| = \|c-y\| = \frac{d}{2}$. If $\tilde z \neq z$, by the strict convexity of $E$
\[
\|c-x\| = \|\frac{1}{2}(z + \tilde z) - x\| < \frac{d}{2}
\]
since $\|z-x\| = \|\tilde z -x\| = \frac{d}{2}$. This contradiction asserts that $\tilde z =z$. So we have shown that any weak cluster point of $\{(T_t)^n(z)\}$ is equal to $z$. Hence
\[
\text{wk-}\lim_{n\to\infty} (T_t)^n (z) =z.
\]
By the weak continuity of $T_t$ we end up with $T_t z =z$ for each $t\in S$, or $z\in F(S) = C_0$. Therefore $C_0$ is convex. The proof is complete.

\end{proof}

 Let $C$ be any nonempty subset of a Banach space $E$ and $\Sc$ be a representation of $S$ on $C$. We are interested in when $A_C(\Sc) \neq \emptyset$. For this purpose it is reasonable to assume that $\{T_sx: s\in S\}$ is bounded for all $x\in C$, since otherwise $A_C(\Sc)$ must be empty. If the representation is nonexpansive, then it is readily seen that $\{T_sx: s\in S\}$ is bounded for all $x\in C$ if there is a point $c\in C$ such that $\{T_sc: s\in S\}$ is bounded. When $E$ is reflexive this condition implies further that $\{T_sx: s\in S\}$ is weakly precompact in $E$ for each $x\in C$. 
 
 \section{nonexpansive mappings on Hilbert spaces}\label{Hilbert}

Let $H$ be a Hilbert space over $\Cbb$. The inner product of $x,y\in H$ is denoted by $\la x \,|\, y \ra$. The following are elementary properties of a Hilbert space.
\begin{equation}\label{identity 1}
\text{Re}\la x+y \,|\, x-y \ra = \|x\|^2 - \|y\|^2 \quad (x,y\in H);
\end{equation}
\begin{equation}\label{identity 2}
\|\lamd x +(1-\lamd) y\|^2 = \lamd \|x\|^2 + (1-\lamd) \|y\|^2 - \lamd(1-\lamd)\|x-y\|^2  
\end{equation}
for $x,y\in H$, $\lamd\in \Rb$.
 
 For a semigroup representation in a Hilbert space  the following holds.

\begin{lemma}\label{As Hilbert}
Let $C$ be a nonempty subset of a Hilbert space $H$ and $\Sc$ be a representation of $S$ on $C$ as nonexpansive self mappings. Suppose that $\{T_sc: s\in S\}$ is bounded for some $c\in C$. Let $X$ be a closed subspace of $\ell^\infty(S)$ containing the constant functions and being invariant under left translations.  If $X$ has a left invariant mean $\mu$ and if $y_c\in X$ for each $y\in H$, where 
\[
y_c(s) = \la T_sc \,|\, y\ra \quad s\in S,
\]
 then $T_\mu c \in A(\Sc)$ and so $A(\Sc) \neq \emptyset$.
\end{lemma}


\begin{proof}
From the discussion in the end of Section~\ref{preliminaries}, $a =T_\mu c\in H$ is well-defined. We show $\|a - T_t x\| \leq \|a - x\|$ for all $x\in C$. In fact, using identity (\ref{identity 1}) we have
\begin{align*}
\|a - T_t x\|^2 - \|a - x\|^2 &= \text{Re}\la 2a -T_tx -x \p x- T_tx\ra \\
                              &= \text{Re}\mu_s(\la 2T_s c -T_tx -x \p x- T_tx\ra)\\
                              &= \mu_s(\text{Re}\la 2T_sc -T_tx -x \p x- T_tx\ra) \\
                              &= \mu_s(\|T_{s}c - T_tx\|^2 -\|T_s c - x\|^2)\\
                              &= \mu_s(\|T_{ts}c - T_tx\|^2) -\mu_s(\|T_s c - x\|^2)\\
                              & \leq \mu_s(\|T_sc -x\|^2) -\mu_s(\|T_s c - x\|^2) =0.
\end{align*}
Note that the last inequality holds because $\mu$ is left invariant. Therefore $\|a - T_t x\|^2 \leq \|a - x\|^2$ for all $t\in S$, and hence $a\in A(\Sc)$.

\end{proof}

In particular, let $S$ be a semitopological semigroup, then for $X$ being $\ell^\infty(S)$, $C_b(S)$, $LUC(S)$, $AP(S)$ or $WAP(S)$ we derive the following result.

\begin{cor}\label{nonempty As}
Let $C$ be a nonempty subset of a Hilbert space $H$ and $\Sc$ be a representation of a semitopological semigroup $S$ on $C$ as nonexpansive self mappings. Suppose that $\{T_sc: s\in S\}$ is bounded for some $c\in C$. Then $A(\Sc) \neq \emptyset$ if any of the following conditions holds.
\begin{enumerate}
\item $C_b(S)$ has a left invariant mean and the mapping $s\mapsto T_sc$ is continuous from $S$ into $(C,\text{wk})$;
\item $S$ is left amenable and the action of $S$ on $C$ is weakly jointly continuous;
\item $AP(S)$ has a left invariant mean and the action of $S$ on $C$ is weakly separately continuous and weakly equicontinuous continuous;
\item $WAP(S)$ has a left invariant mean and the action of $S$ on $C$ is weakly separately continuous and weakly quasi-equicontinuous.
\end{enumerate}
\end{cor}

\begin{proof}
By Lemma~\ref{As Hilbert} it suffices to verify that for $X$ being 
 $C_b(S)$, $LUC(S)$, $AP(S)$ or $WAP(S)$ under the corresponding condition of (1)--(4), we have $y_c \in X$ for all $y\in H$. This is clear for case (1) 
  In other cases we may assume $C = \overline{\{T_sc: s\in S\}}^\text{wk}$, the closure of $\{T_sc: s\in S\}$ in the weak topology of $H$. Then $C$ is compact in the weak topology of $H$. In case (2), $LUC(S)$ has a left invariant mean. Routine computation shows that $y_c \in LUC(S)$ if the representation of $S$ on $C$ is weakly jointly continuous. In case (3) $y_c\in AP(S)$ is due to \cite[Lemma 3.1]{Lau73}. In case (4) $y_c\in WAP(S)$ is due to \cite[Lemma 3.2]{L-Z1}.
\end{proof}


Combining Lemma~\ref{As to Fs} and Corollary~\ref{nonempty As}, we obtain a common fixed point theorem for an acting on an unbounded convex set of a Hilbert space.

\begin{thm}\label{fpp Hilbert}
 Let $C$ be a nonempty closed convex subset of a Hilbert space $H$ and $\Sc$ be a representation of $S$ on $C$ as nonexpansive self mappings. Suppose that $\{T_sc: s\in S\}$ is bounded for some $c\in C$. Then $F(\Sc) \neq \emptyset$ if any of the conditions (1)-(4) in Lemma~\ref{nonempty As} holds
\end{thm}

It is well known that if a semitopological semigroup $S$ is left reversible then $AP(S)$ has a left invariant mean \cite{Lau73}. For a discrete semigroup $S$ the left reversibility also implies that $WAP(S)$ has a left invariant mean \cite{Hsu, L-Z1}. We are interested in whether $y_c\in AP(S)$ or $\in WAP(S)$ for $y\in H$ and $c\in C$ if $S$ is a left reversible semitopological semigroup acting on a set $C$ of a Hilbert space $H$.
For an affine representation we have the following.

\begin{prop}\label{affine}
Let $C$ be a nonempty convex subset of a Hilbert space $H$ and $\Sc$ be a representation of a semitopological semigroup $S$ on $C$ as nonexpansive affine self mappings such that the mapping $s\mapsto T_s x$ from $S$ into $(C,\text{wk})$ is continuous for each $x\in C$. Suppose that $\{T_sc: s\in S\}$ is bounded for some $c\in C$. If $WAP(S)$ has a left invariant mean $\mu$, then 
 $F(\Sc)\neq \emptyset$
\end{prop}

\begin{proof}
By Lemma~\ref{As to Fs} and Lemma~\ref{As Hilbert} we only need to show $y_c\in WAP(S)$ for all $y\in H$. For this purpose we may assume $C=\overline{\text{co}}^\text{wk}(\{T_sc: s\in S\})$, which is weakly compact (note that a norm continuous affine mapping is always weakly continuous). We shall show that the representation of $S$ on $C$ is quasi-equicontinuous when $C$ is equipped with the weak topology of $H$. 

Suppose that $(s_i)\subset S$ is a net satisfying $T_{s_i}(x) \overset{\text{wk}}{\to} T(x)$ for each $x\in C$. We show $T$ is weak-weak continuous. If this were not true, then there would be a net $(x_j)\subset C$ such that $x_j\overset{\text{wk}}{\to} x$ but $T(x_j)\overset{\text{wk}}{\nrightarrow} T(x)$. Then there would exist $z\in H$, $\ep>0$ and a subnet of $(x_j)$, still denoted by $(x_j)$, such that $\|z\| =1$ and 
\[
\text{Re} [\la T(x_j) -T(x) \,|\, z\ra] > \ep
\]
for all $j$. By Mazur's Theorem, there is a net $(x_\lambda)\subset \text{co}(x_j)$ such that $x_\lambda \to x$ in norm. We certainly still  have 
\[
\text{Re} [\la T(x_\lambda) -T(x) \,|\, z\ra] > \ep
\]
for all $\lambda$. On the other hand
\[
\la T(x_j) -T(x) \,|\, z\ra = \lim_i \la T_{s_i}(x_j) - T_{s_i}(x)\, |\, z\ra.
\]
By the nonexpansiveness of $T_{s_i}$ we would have
\[
\|T(x_j) - T(x)\| \leq \|x_j-x\| \to 0.
\]
This contradiction shows that $T$ is weak-weak continuous. Thus the representation of $S$ on $C$ is weakly quasi-equicontinuous. So $y_c\in WAP(S)$ for all $y\in H$ from \cite[Lemma 3.2]{L-Z1}. The proof is complete.
 
\end{proof}

\begin{cor}
If $S$ is a discrete left reversible semigroup, then, 
 a nonexpansive affine representation of 
 $S$ on a convex set $C$ of a Hilbert space has 
 a common fixed point if $\{T_sc: s\in S\}$ is bounded for some $c\in C$.
\end{cor}

We now consider when $y_c\in AP(S)$ assuming $S$ is left reversible.

\begin{lemma}\label{TB=B}(\cite[Lemma~3.4]{L-Z2})
Let $S$ be a left reversible semitopological semigroup and $\Sc =\{T_s: s\in S\}$ a representation of $S$ as jointly continuous self mappings on a compact Hausdorff space $K$. Then there is a nonempty compact subset $B$ of $K$ such that $T_s(B) = B$ for all $s\in S$.
\end{lemma}

\begin{lemma}\label{norm cpt}
Let $S$ be a separable semitopological semigroup that acts on a weakly compact subset $K$ of a Banach space $E$ as weakly separately continuous and norm nonexpansive mappings. Suppose that $F$ is a minimal nonempty weakly compact $S$-invariant subset of $K$ satisfying $sF = F$ for all $s\in S$. Then $F$ is norm compact.
\end{lemma} 
\begin{proof}
This follows from the same proof of (\cite[Lemma~3.3]{L-Z1}), where the convex assumption for $K$ was not used in the argument and so is removable.

 
\end{proof}

\begin{thm}\label{left rev}
Let $S$ be a left reversible and separable semitopological semigroup, and let $\Sc=\{T_s: s\in S\}$ be a representation of $S$ on a weakly closed subset $C$ of a Hilbert space $H$ as norm nonexpansive and weakly jointly continuous self mappings. If there is $c\in C$ such that $\{T_s c: s\in S\}$ is bounded, then $A_C(\Sc) \neq \emptyset$. In particular, $F(\Sc) \neq \emptyset$.
\end{thm}
\begin{proof}
Since $S$ is left reversible, by Lemma~\ref{TB=B} and Zorn's Lemma there is a minimal weakly compact $S$-invariant subset $F$ of $\overline{\{T_s c: s\in S\}}^w\subset C$ such that $T_s(F) = F$ for all $s\in S$. From Lemma~\ref{norm cpt} $F$ is actually norm compact. Then it follows from \cite[Lemma~3.1]{Lau73} the function $y_c$ is in $AP(S)$ for $c\in F$ and $y\in H$. On the other hand $AP(S)$ has a left invariant mean if $S$ is left invariant. Therefore, by Lemma~\ref{As Hilbert} $A_C(\Sc) \neq \emptyset$. Then applying Lemma~\ref{As to Fs} we have $F(\Sc)\neq \emptyset$.

\end{proof}

\begin{lemma}\label{convex}
Let $C$ be a nonempty subset of a Hilbert space $H$ and $\Sc$ be a representation of $S$ on $C$. If $A_C(\Sc) \neq \emptyset$ then $A_C(\Sc)$ is closed and convex.
\end{lemma}

\begin{proof}
From the definition it is evident that $A_C(\Sc)$ is closed. Let $a,b\in A_C(\Sc)$ and $0\leq \lamd \leq 1$. Then for $c\in C$ and $t\in S$
\begin{align*}
\|\lamd a +(1-\lamd) b -T_t c\|^2 & = \lamd\|a-T_t c\|^2 + (1-\lamd)\|b-T_t c\|^2 - \lamd(1-\lamd)\|a-b\|^2\\
                                  & \leq \lamd\|a- c\|^2 + (1-\lamd)\|b- c\|^2 - \lamd(1-\lamd)\|a-b\|^2\\
                                  & = \|\lamd a +(1-\lamd) b - c\|^2.
\end{align*}
Therefore $\lamd a +(1-\lamd) b \in A_C(\Sc)$.

\end{proof}

Suppose that $S$ is right reversible, i.e. $\overline{Sa}\cap \overline{Sb} \neq \emptyset$. We may define a partial order on $S$ by letting $a \leq b$ if $\overline{Sa}\cup\{a\} \supseteq \overline{Sb}\cup\{b\}$. Let $T$ be a mapping from $S$ into a topological space $X$. We say $\{T(s)\}$ converges to $y\in X$ if the limit $\lim_sT(s) = y$ holds when $s$ increases in this order.

\begin{lemma}\label{CD}
Let $C$ and $D$ be  non-empty sets in a Hilbert space $H$ with $D$ closed and convex. Let $S$ be a right reversible semitopological semigroup and $\Sc=\{T_s: s\in S\}$ a representation of $S$ on $C$. Suppose that $x$ is a point in $C$ such that the mapping $s\mapsto T_s x$ is continuous from $S$ into $C$ and
\[
\|T_{ts}x - z\| \leq \|T_sx-z\| \quad (z\in D, s,t\in S).
\]
Then $\{P(T_sx)\}$ converges strongly to some $z_0\in D$, where $P$: $H\to D$ is the metric projection. 
\end{lemma}

\begin{proof}
We first notice that by the definition of $P$
\[
\|y-P(y)\| \leq \|y-z\| \quad (y\in H, z\in D).
\]
If $a,b\in S$, $a\geq b$, then either $a=b$ or there is a net $(s_i)\subset S$ such that $a =\lim_i s_ib$. In the latter case,
\[
\|T_a x - z\|=\lim_i \|T_{s_ib}x -z\| \leq \|T_bx-z\| \quad (z\in D)
\]
by the assumption. In particular, $\|T_ax - P(T_bx)\|\leq \|T_bx - P(T_bx)\|$ and hence
\[
\|T_ax - P(T_ax)\|\leq \|T_bx - P(T_bx)\| \quad (a\geq b).
\]
 So the numerical net $\{\|T_sx-P(T_sx)\|\}$ is decreasing as $s$ increases. Thus the limit $\lim_s \|T_s x -P(T_sx)\|$ exists.
 
On the other hand, for $y\in H$ and $z\in D$
\[
\|y-z\|^2 = \|y-P(y)\|^2 + \|P(y)-z\|^2 + 2\text{Re}\la y-P(y) \,|\, P(y)-z \ra.
\]
 Using (\ref{mp}) we immediately get
\[
\|P(y) -z\|^2 \leq \|y - z\|^2 - \|P(y)-y\|^2 \quad (y\in H, z\in D.)
\]
Therefore,
\begin{align*}
\|P(T_ax) - P(T_bx)\|^2 &\leq \|T_ax - P(T_bx)\|^2 - \|P(T_ax) - T_ax\|^2\\
                          &\leq \|T_bx -P(T_bx)\|^2 - \|T_ax-P(T_ax)\|^2 .
\end{align*}
This shows that $\{P(T_sx)\}$ is a Cauchy net in $H$. It then converges to some $z_0 \in D$.

\end{proof}

We note that Lemma~\ref{CD} was proved in \cite[Proposition~2.4]{Lau85} for the case $D=F(S)$ when the representation is nonexpansive. 

 We now use Lemma~\ref{CD} to prove an ergodic theorem for representations of a semitopological semigroup in a Hilbert space. 

\begin{thm}\label{ergodic}
Let $C$ be a nonempty closed subset of a Hilbert space $H$, let $\Sc=\{T_s: s\in S\}$ be a representation of a right reversible semitopological semigroup $S$ on $C$ as separately continuous nonexpansive mappings, and let $X$ be  the left invariant subspace of $\ell^\infty(S)$ in any of the cases listed below. Suppose that $X$ has a LIM $\mu$ and suppose that $\Sc c =\{T_sc: s\in S\}$ is bounded for some $c\in C$. Then $A_C(\Sc)$ is a nonempty closed convex subset of $H$, $T_\mu x\in A_C(\Sc)$ for all $x\in C$. In particular, $F(\Sc)\neq \emptyset$. If in addition $\mu$ is also a RIM on $X$, then $\lim_t P_{A_C(\Sc)}(T_t x) = T_\mu x$.
\begin{enumerate}
\item $X=C_b(S)$ and the mapping $s\mapsto T_sc$ is continuous from $S$ into $(C,\text{wk})$;
\item $X=LUC(S)$, and the action of $S$ on $C$ is weakly jointly continuous;
\item $X=AP(S)$  and the action of $S$ on $C$ is weakly equicontinuous continuous;
\item $X=WAP(S)$ and the action of $S$ on $C$ is weakly quasi-equicontinuous;
\item $S$ is left reversible and separable and the action of $S$ on $C$ is weakly 
equicontinuous continuous (for this case $\mu$ is taken to be a LIM for $AP(S)$)
\end{enumerate}
\end{thm}

\begin{proof}
First, by the assumption the mapping $s\mapsto T_s x$ is continuous in norm, hence is also continuous in the weak topology for $C$. So in the cases (3) and (4) the action of $S$ on $C$ is automatically weakly separately continuous. 
 From Lemma~\ref{nonempty As} Theorem~\ref{left rev} and Lemma~\ref{convex}, we know that, in each case listed, $A_C(\Sc) $ is nonempty and is closed convex. In particular, $F(\Sc) \neq \emptyset$. We  also note that from the nonexpansiveness of the representation, $\Sc x$ is bounded for all $x\in C$. 
 
 We now show that the limit $\lim_t P_{A_C(\Sc)}(T_t x) = T_\mu x$ holds if $\mu$ is an invariant mean on $X$.

 By Lemma~\ref{CD} with $D = A_C(\Sc)$, $\lim_t P_{A_C(\Sc)}(T_t x)$ exists. Assume the limit is $u\in A_C(\Sc)$. By the property of $P_{A_C(\Sc)}$ we have
\[
\text{Re}\left[\la T_t x -  P_{A_C(\Sc)}(T_t x) \, |\,  P_{A_C(\Sc)}(T_t x) -y \ra \right] \geq 0 \quad (t\in S, y\in A_C(\Sc)).
\]
This leads to 
\[
\text{Re}\left[\la T_t x -  P_{A_C(\Sc)}(T_t x) \, |\,  y - u \ra \right] \leq \text{Re}\left[\la T_t x -  P_{A_C(\Sc)}(T_t x) \, |\,  P_{A_C(\Sc)}(T_t x) -u \ra \right] .
\]
Since $\{T_s x: s\in S\}$ is bounded, it is evident that  $\{P_{A_C(\Sc)}(T_s x): s\in S\}$ is bounded too. We assume 
\[
\|T_s x\| + \|P_{A_C(\Sc)}(T_s x)\| \leq M.
\]
 Then
\[
\text{Re}\left[\la T_t x -  P_{A_C(\Sc)}(T_t x) \, |\,  y - u \ra\right]  \leq M \|P_{A_C(\Sc)}(T_t x) -u\| \quad (t\in S, y\in A_C(\Sc)) .
\]
Apply $\mu$ on both sides. We have
\[
\text{Re}\left[\mu_s\left(\la T_{st} x -  P_{A_C(\Sc)}(T_{st} x) \, |\,  y - u \ra\right) \right] \leq M\mu_s\left( \|P_{A_C(\Sc)}(T_{st} x) -u\|\right)
\]
for $t\in S$ and $y\in A_C(\Sc)$.
By the right invariance of $\mu$ we get
\[
\text{Re}\left[\mu_s\left(\la T_{\mu} x -  P_{A_C(\Sc)}(T_{st} x) \, |\,  y - u \ra\right)\right]  \leq M\mu_s\left( \|P_{A_C(\Sc)}(T_{st} x) -u\|\right) 
\]
for $t\in S$ and  $y\in A_C(\Sc)$.
Since $\lim_t P_{A_C(\Sc)}(T_{st} x)$ converges to $u$ uniformly in $s$ by the definition of the order on $S$, we finally get
\[
\text{Re}\left[\la T_\mu x -u\, |\, y-u\ra\right]  \leq 0 \quad (y\in A_C(\Sc)). 
\]
In particular this is true for $y = T_\mu x$. Thus we derive $\|T_\mu x - u\|^2 \leq 0$ or $T_\mu x = u$.

\end{proof}

\begin{remark}
A special case of our Theorem~\ref{ergodic} for a discrete commutative semigroup representation was obtained in \cite{A-T} (see Theorem~4.1 there). 
\end{remark}

\begin{remark} We indicate here several important cases regarding the existence of a LIM.
\begin{enumerate}
\item
When $S$ is discrete semigroup, the following implication diagram is known (see \cite{L-Z1}).
\begin{equation*}
\begin{matrix}
S \text{ left amenable}\\
\Downarrow \quad \mathrel{\not \Uparrow}\\
S \text{ left reversible}\\
\Downarrow \quad \mathrel{\not \Uparrow}\\
WAP(S) \text{ has LIM}\\
\Downarrow \quad \mathrel{\not \Uparrow}\\
AP(S) \text{ has LIM}\\
\end{matrix}
\end{equation*}
\item If $G$ is a topological group, then $WAP(G)$ always has a left invariant mean by the Ryll-Nardzewskii fixed point theorem.
\item If $G$ is a locally compact group, then $C_b(G)$ has a LIM if and only if $G$ is left amenable.
\end{enumerate}
\end{remark}

It is well known that $X$ has a left (right, or two-sided) invariant mean if and only if there is a net $(\mu_\al)$ of means on $X$ such that $\mu_\al - \ell^*_s\mu_\al \to 0$ (resp. $\mu_\al - r^*_s\mu_\al \to 0$, or both limits hold) in the weak* topology of $X^*$ for all $s\in S$. Moreover, if $(\mu_\al)$ is such a net then a subnet of it will converge weak* to a left (resp. right, or two-sided) invariant mean $\mu$. On the other hand, if $\mu_\al \overset{\text{wk*}}{\to} \mu$ then $T_{\mu_\al} x\overset{\text{wk}}{\to} T_\mu x$. 

\begin{cor}
Under the condition of Theorem~\ref{ergodic} if $(\mu_\al)$ is a net of means on $X$ such that 
\[
\mu_\al - \ell^*_s\mu_\al \overset{\text{wk*}}{\to} 0 \text{ and } \mu_\al - r^*_s\mu_\al \overset{\text{wk*}}{\to} 0 \quad  (s\in S),
\]
then for each $x\in C$, $\{T_{\mu_\al}(x)\}$ converges weakly to 
\[
u_x: = \lim_t P_{A_C(\Sc)}(T_t x)\in A_C(\Sc).
\]
Moreover, if $\mu_\al \to \mu$ in norm, then $(T_{\mu_\al}x)$
also converges to $u_x$ in norm.
\end{cor}

\begin{proof}
First assume $\mu_\al \overset{\text{wk*}}{\to} \mu$. Then $\mu$ is an invariant mean on $X$.
Since $T_\mu x= u_x = \lim_t P_{A_C(\Sc)}(T_t x)$ for $x\in C$ from Theorem~\ref{ergodic}, and since $T_{\mu_\al} x\overset{\text{wk}}{\to} T_\mu x$, we conclude immediately that $\{T_{\mu_\al} x\}$ converges weakly to $u_x$ for all $x\in C$.

In general, the above shows that for every weak* convergent subnet $(\mu_\beta)$ of $(\mu_\al)$
\[
T_{\mu_\beta} x\overset{\text{wk}}{\to} u_x = \lim_t P_{A_C(\Sc)}(T_t x).
\]
This implies $T_{\mu_\al} x\overset{\text{wk}}{\to} u_x$.

Moreover, if $\mu_\al \to \mu$ in norm, then $(T_{\mu_\al}x)$
also converges to $u_x$ in norm. In fact,
\[
\|T_{\mu_\al}x - T_\mu x\| = \sup_{f\in H_1}|(\mu_\al -\mu)(\la T_sx\,|\, f \ra| \leq M\|\mu_\al - \mu\|
\]
where $M > 0$ is a constant satisfying $\|T_s x\| \leq M$ for all $s\in S$. We have
\[
\|T_{\mu_\al}x - u_x\| = \|T_{\mu_\al}x - T_\mu x\| \overset{\al}{\to} 0.
\]
The proof is complete.

\end{proof}

\begin{cor}
Let $S$ be a reversible discrete semigroup (i.e. $S$ be both left and right reversible). If $\Sc = \{T_s: s\in S\}$ is a representation of $S$ on a nonempty closed subset $C$ of a Hilbert space $H$ such that the action of $S$ on $C$ is weakly quasi-continuous, then $F(\Sc)\neq \emptyset$.
\end{cor}
\begin{proof}
This follows from Theorem~3.11 and a result of R. Hsu in \cite{Hsu}, where he showed that when $S$ is left reversible and discrete, then $WAP(S)$ has a LIM.
\end{proof}

\section{generalized hybrid mappings in Hilbert spaces}\label{hybrid}

In this section we aim to extend the results in the previous section to semigroups of mappings generated by so called generalized hybrid mappings. Let $E$ be a Banach space and $C\subset E$. We call a mapping $T$: $C\to C$ a \emph{generalized hybrid mapping} \cite{K-T-Y} if there are numbers $\al, \beta \in \Rb$ such that
\[
\al\|Tx - Ty\|^2 + (1-\al)\|x - Ty\|^2 \leq \beta \|Tx - y\|^2 + (1 - \beta)\|x - y\|^2.
\]
for all $x,y \in C$.

When $(\al,\beta) = (1,0)$, this indeed defines a nonexpansive mapping. However, in general the composite of two generalized hybrid mappings is usually no long generalized hybrid.  Also a generalized hybrid mapping may be discontinuous.

\begin{lemma}\label{As hybrid}
Let $C$ be a nonempty subset of a Hilbert space $H$ and $\Sc=\{T_s: s\in S\}$ be a representation of $S$ on $C$. Suppose that $S$ is generated by a subset $\Lambda$ and $T_s$ is a generalized hybrid mapping on $C$ for each $s\in \Lambda$ and suppose that $\{T_sc: s\in S\}$ is bounded for some $c\in C$. Let $X$ be a closed subspace of $\ell^\infty(S)$ containing the constant functions and being invariant under left translations.  If $X$ has a left invariant mean $\mu$ and if $y_c\in X$ for each $y\in H$, where 
\[
y_c(s) = \la T_sc \,|\, y\ra \quad s\in S,
\]
 then $T_\mu c \in A(\Sc)$ and so 
 $F(\Sc) \neq \emptyset$.
\end{lemma}

\begin{proof}
Still from the discussion in the end of Section~\ref{preliminaries}, $a =T_\mu c\in H$ is well-defined. Following the proof of Lemma~\ref{As Hilbert} we have
\begin{align*}
\|a - T_t x\|^2 - \|a - x\|^2 &= \mu_s(\|T_{s}c - T_tx\|^2 -\|T_s c - x\|^2)\\
                              &= \mu_s(\|T_{s}c - T_tx\|^2) -\mu_s\|T_s c - x\|^2) \quad (t\in S, x\in C).
\end{align*}

Now let $t\in \Lambda$. We have
\begin{align*}
\mu_s(\|T_{s}c - T_tx\|^2) &= \mu_s(\al \|T_{s}c - T_tx\|^2) + (1-\al)\|T_{s}c - T_tx\|^2)\\
                           &= \mu_s(\al \|T_{ts}c - T_tx\|^2) + (1-\al)\|T_{s}c - T_tx\|^2)\\
                           & \leq \mu_s(\beta \|T_{ts}c - x\|^2) + (1-\beta)\|T_{s}c - x\|^2)\\
                           & = \mu_s(\beta \|T_{s}c - x\|^2) + (1-\beta)\|T_{s}c - x\|^2)\\
                           & = \mu_s(\|T_{s}c - x\|^2) \quad(x\in C)
\end{align*}
for some $\al,\beta \in \Rb$. Thus $\|a - T_t x\|^2 - \|a-x\|^2 \leq 0$ for $t\in \Lambda$ and $x\in C$. Since $\Lambda$ generates $S$, the inequality still holds for all $t\in S$. Therefore, $a= T_\mu c$ is an attractive point of $S$ for $C$.

\end{proof}

With the above lemma we may establish analogues of Theorems~\ref{fpp Hilbert} and \ref{ergodic} for semigroups generated by generalized hybrid mappings with the same proofs.

\begin{thm}\label{hybrid As}
Let $C$ be a nonempty subset of a Hilbert space $H$ and $\Sc$ be a representation of a semitopological semigroup $S$ on $C$.  Suppose that $S$ is generated by a subset $\Lambda$ and $T_s$ is a generalized hybrid mapping on $C$ for each $s\in \Lambda$, and suppose that $\{T_sc: s\in S\}$ is bounded for some $c\in C$. Then 
$F(\Sc) \neq \emptyset$ if any of the following conditions holds.
\begin{enumerate}
\item $C_b(S)$ has a left invariant mean and the mapping $s\mapsto T_sc$ is continuous from $S$ into $(C,\text{wk})$;
\item $S$ is left amenable and the action of $S$ on $C$ is weakly jointly continuous;
\item $AP(S)$ has a left invariant mean and the action of $S$ on $C$ is weakly separately continuous and weakly equicontinuous continuous;
\item $WAP(S)$ has a left invariant mean and the action of $S$ on $C$ is weakly separately continuous and weakly quasi-equicontinuous.
\end{enumerate}
\end{thm}

\begin{lemma}\label{ergodic hybrid}
Let $C$ be a nonempty closed subset of a Hilbert space $H$, let $\Sc=\{T_s: s\in S\}$ be a representation of a right reversible semitopological semigroup $S$ on $C$ such that the mapping $s\mapsto T_s x$ is continuous from $S$ into $C$ with the norm topology for each $x\in C$. Suppose that $S$ is generated by a subset $\Lambda$ and $T_s$ is a generalized hybrid mapping on $C$ for each $s\in \Lambda$, and suppose that $\{T_sc: s\in S\}$ is bounded for some $c\in C$. Let $X$ be  the left invariant subspace of $\ell^\infty(S)$ in any of the cases listed below such that $X$ has a two-sided invariant mean $\mu$. Then $A_C(\Sc)$ is a nonempty closed convex subset of $H$, $T_\mu x\in A_C(\Sc)$ for all $x\in C$ and $\lim_t P_{A_C(\Sc)}(T_t x) = T_\mu x$.
\begin{enumerate}
\item $X=C_b(S)$ and the mapping $s\mapsto T_sc$ is continuous from $S$ into $(C,\text{wk})$;
\item $X=LUC(S)$, and the action of $S$ on $C$ is weakly jointly continuous;
\item $X=AP(S)$  and the action of $S$ on $C$ is weakly equicontinuous continuous;
\item $X=WAP(S)$ and the action of $S$ on $C$ is weakly quasi-equicontinuous;
\end{enumerate}
\end{lemma}

\begin{proof}
Due to Theorem~\ref{hybrid As} $A_C(\Sc) \neq \emptyset$. This in turn implies that $\{T_sx: s\in S\}$ is bounded for each $x\in C$. Then using the argument for the proof of Theorem~\ref{ergodic} we get the result. 
\end{proof}

Applying Lemma~\ref{As to Fs} we then obtain a fix point theorem for generalized hybrid mappings in Hilbert spaces as follows.

\begin{thm}
Suppose that the condition of Lemma~\ref{ergodic hybrid} holds with $C$ being convex. Then there is a common fixed point for $S$ in $C$.
\end{thm}

\section{Some remarks and open problems}\label{problems}

A semitopological semigroup $S$ is \emph{extremely left amenable} if $LUC(S)$ has a multiplicative left invariant mean. If $S$ is a locally compact group, then $S$ is extremely left amenable only when $S$ is a singleton \cite{G-L70}. However, a non-trivial topological group which is not locally compact can be extremely left amenable. In fact, let $S$ be the group of unitary operators on an infinite dimensional Hilbert space with the strong operator topology, then $S$ is extremely left amenable \cite{G-M}; In \cite[Theorem~3.2]{L-Z3} the authors showed that an F-algebra $A$ is left amenable if and only if the semigroup of normal positive functions of norm 1 on $A^*$ is extremely left amenable. For more examples we refer to \cite{Pestov, Lau-Ludwig}.

\begin{problem}
Suppose that $S$ is extremely left amenable and $C$ is a weakly closed subset of a Banach space $E$, and suppose that $\Sc$ is a weakly continuous and norm nonexpansive representation of $S$ on $C$ such that $Sc = \{T_sc: s\in S\}$ is relatively weakly compact for some $c\in C$. Does $C$ contain a fixed point for $S$? \label{Prob extreme}
\end{problem}
We know that the answer is ``yes'' when $S$ is discrete. Indeed, in this case, for each finite subset $\sigma$ of $S$ there is $s_\sigma \in S$ such that $ss_\sigma =s_\sigma$ for all $s\in \sigma$ by a theorem of Granirer's \cite{Granirer} (see also \cite[Theorem~4.2]{L-Z2} for a short proof). Consider the net $\{s_\sigma c\}$. By the relative weak compactness of $Sc$, there is $z\in \overline{Sc}^\text{ wk}\subset C$ such that (replaced by a subnet if necessary) wk-$\lim_{\sigma}s_\sigma c = z$. Then, as readily checked, $T_sz = z$ for all $s\in S$ by the weak continuity of the $S$ action on $C$.

More generally, the answer to Problem~\ref{Prob extreme} is still affirmative (even without the norm nonexpansiveness assumption) if the representation is jointly continuous when $C$ is equipped with the weak topology of $E$. This is indeed a consequence of \cite[Theorem~1]{Mitch_LUC} or \cite[Theorem 5.4(a)]{L-Z2}.

\begin{problem}
 Let $C$ be a nonempty closed convex subset of the sequence space $c_0$ and $\Sc$ be a representation of a commutative semigroup $S$ as nonexpansive mappings on $C$. Suppose that $\{T_sc: s\in S\}$ is relatively weakly compact for some $c\in C$. Is $F(\Sc) \neq \emptyset$?\label{Prob c0}
\end{problem}
One may not drop the weak compactness condition on the orbit of $c$. T. D. Benavides has shown that for any unbounded subset $C$ of $c_0$ there is an nonexpansive mapping $T$ on $C$ such that $C$ has no fixed point for $T$ \cite{Benavides}. In fact, even $C$ is bounded, without weak compactness of an orbit the answer to Problem~\ref{Prob c0} will still be negative. For example, on the unit ball of $c_0$ define $T((x_i)) = (1,x_1, x_2, \cdots)$. Then $T$ is nonexpansive, and obviously $T$ has no fixed point in the unit ball.

\begin{problem}
Does Lemma~\ref{As to Fs} still hold when $E$ is merely a strictly convex Banach space?
\end{problem}

\begin{problem}
In any of the cases studied in Theorems~\ref{fpp Hilbert}, \ref{affine}, \ref{left rev} or \ref{ergodic}, does the converse hold?
\end{problem}

The authors are grateful to the referee for drawing their attention to the article \cite{T-T}.

\end{document}